\documentclass[11pt,a4paper]{amsart}
\usepackage{amsmath}

\usepackage{amssymb,latexsym}

\newtheorem{theorem}{Theorem} [section]
\newtheorem{lemma}{Lemma} [section]
\newtheorem{proposition}{Proposition} [section]
\newtheorem{corollary}{Corollary} [section]
\newtheorem{definition}{Definition} [section]
 
\newtheorem{remark}{Remark}[section]

\begin{document}

\title[Critical potentials for Schr\"{o}dinger operators]{Critical potentials of the eigenvalues and eigenvalue gaps
of Schr\"{o}dinger operators}

\author{ Ahmad El Soufi   and   Nazih Moukadem }

\address{A. El Soufi : Laboratoire de Mathématiques et Physique Théorique, 
UMR CNRS 6083, Université de Tours, Parc de Grandmont, F-37200
Tours France}
\email{elsoufi@univ-tours.fr}
\address{N. Moukadem : Département de Mathématiques, Université Libanaise, 
Faculté des Sciences III, Tripoli Liban}

\keywords{eigenvalues, Schrödinger operator, extremal potential, extremal gap}
\subjclass{35J10, 35P15, 49R50, 58J50 }

\begin{abstract}
Let $M$ be a compact Riemannian manifold with or without boundary, and let $-\Delta $ be its Laplace-Beltrami operator. For any bounded scalar potential $q$, we denote by $\lambda_i(q)$ the $i$-th eigenvalue of the Schr\"{o}dinger type operator $-\Delta + q$ acting on functions with Dirichlet or Neumann boundary conditions in case $\partial M \neq \emptyset$. We investigate critical potentials of the eigenvalues $\lambda_i$ and the eigenvalue gaps $G_{ij}=\lambda_j -\lambda_i$ considered as functionals on the set of bounded potentials having a given mean value on $M$. We give necessary and sufficient conditions for a potential $q$ to be critical or to be a local minimizer or a local maximizer of these functionals. For instance, we prove that a potential $q \in L^\infty (M)$ is critical for the functional $\lambda_2$ if and only if, $q$ is smooth, $\lambda_2( q)=\lambda_3( q)$ and there exist second eigenfunctions $f_1 ,\ldots,f_k$ of $-\Delta + q$ such that $\Sigma_j f_j^2 = 1$. In particular, $\lambda_2$ (as well as any $\lambda_i$) admits no critical potentials under Dirichlet Boundary conditions. 
Moreover, the functional $\lambda_2$ never admits locally minimizing potentials.

\end{abstract}

\maketitle

\section {Introduction and Statement of main Results}\label{1}

Let $M$ be a compact connected Riemannian manifold of dimension
$d$, possibly with nonempty boundary $\partial M$, and let $-\Delta $ be its Laplace-Beltrami operator acting on functions with, in the case where $\partial M \neq \emptyset$, Dirichlet or Neumann boundary conditions. In all the sequel, as soon as the Neumann Laplacian will be considered, the boundary of $M$ will be assumed to be sufficiently regular (e.g. $C^1$, but weaker regularity assmptions may suffice, see  \cite{B}) in order to guarantee the compactness of the embedding $H^1(M)   \hookrightarrow  L^2(M)$ and, hence, the compactness of the resolvent of the Neumann Laplacian (note that it is well known, using standard arguments like in  \cite[p.89]{Jo}, that compactness results for Sobolev spaces on Euclidean domains remain valid in the Riemannian setting).

For any bounded real valued potential $q $ on $M$, the Schr\"{o}dinger type operator $-\Delta + q$ has compact resolvent (see \cite[Theorem IV.3.17]{K} and observe that a bounded $q$ leads to a relatively compact operator with respect to $-\Delta $). Therefore, its spectrum consists of a nondecreasing and unbounded sequence of eigenvalues with finite multiplicities:
$$Spec(-\Delta+q)=\{\lambda_1(q)  <  \lambda_2(q) \le \lambda_3(q)\le 
 \cdots \le \lambda_i(  q)\le \cdots \}.$$
Each eigenvalue $ \lambda_i( q)$ can be considered as a (continuous) function of the potential $q\in  L^\infty (M)$ and there are both physical and mathematical motivations to study existence and properties of extremal potentials of the functionals $ \lambda_i$ as well as of the differences, called gaps, between them. A very rich literature is devoted 
to the existence and the determination of maximizing or minimizing potentials for the eigenvalues (especially the fundamental one, $\lambda_1$) and the eigenvalue gaps (especially the first one, $\lambda_2 - \lambda_1$) under various constraints often motivated by physical considerations (see, for instance, \cite{AH, AHS, C, EI, EI1, EK, F, Ha, H, KS, S} and the references therein). Note that, since the function $ \lambda_i$ commutes with constant translations, that is, $ \lambda_i(  q+c)= \lambda_i(  q)+c$, such constraints are necessary. 

Our aim in this paper is to investigate critical points, including "local minimizers" and "local maximizers", of the eigenvalue functionals $q\to  \lambda_i( q)$  and the eigenvalue gap functionals $q\to  \lambda_j(  q)-\lambda_i(  q)$, the potentials $q$ being subjected to the constraint that their mean value (or, equivalently, their integral) over $M$ is fixed.  
All along this paper, the mean value of an integrable function $q$ will be denoted $\bar q$, that is,
$$\bar q ={1\over V(M)} \int_{M} q\, dv,$$
$V(M)$ and $dv$ being respectively the Riemannian volume and the Riemannian volume element of $M$.

Actually, most of the results below can be extended, modulo some slight changes, to the case where this constraint is replaced by the more general one 
$$\int_MF(q)dv = \hbox{constant},$$
where $F:{\mathbb R}\to{\mathbb R}$ is a continuous function such that $F'(x)\neq 0$ if $x \neq 0$, like $F(x)=|x|^\alpha$ or $F(x)=x|x|^{\alpha-1}$ with $\alpha\ge 1$. However, for simplicity and clarity reasons, we preferred to focus only on the mean value constraint. 
Therefore, we fix a constant $c\in {\mathbb R}$ and consider the functionals 
$$ \lambda_i : q \in L_c^\infty (M)\mapsto \lambda_i (q)\in {\mathbb R},$$
where $L_c^\infty (M)=\left\{q\in L^\infty (M)\; \big| \; \bar q = c \right\}$.
The tangent space to  $L_c^\infty (M)$ at any point  $q$ is given by 
$$L_*^\infty (M):=\left\{u \in L^\infty (M) \;\big|  \int_{M} u\, dv = 0\right\}.$$


\subsection{Critical potentials of the eigenvalue functionals} $\quad$

\medskip

\noindent Since it is always nondegenerate, the first eigenvalue gives rise to a differentiable functional in the sense that, for any $q\in L_c^\infty (M)$ and any $u \in L_*^\infty (M)$, the function $t \mapsto \lambda_1 (q+tu)$ is differentiable in $t$. A potential $q \in L_c^\infty (M)$ will be termed \textit{critical} for this functional if ${d\over dt}\lambda_1 (  q +tu)\big|_{t=0} = 0$ for any $u \in L_*^\infty (M)$.

In the case of empty boundary or of Neumann boundary conditions, the constant function 1 belongs to the domain of the operator $-\Delta +q$ and one obtains, as a consequence of the min-max principle, that the constant potential $c$ is a global maximizer of $\lambda_1$ over $L_c^\infty (M)$ (see also \cite{EI} and \cite{H}). 
Constant potential $c$ is actually the only critical one for $\lambda_1$. On the other hand, under Dirichlet boundary conditions, the functional $\lambda_1$ admits no critical potentials in $L_c^\infty (M)$. Indeed, we have the following
\begin{theorem}
\begin{enumerate}
\item Assume that either $\partial M = \emptyset$ or $\partial M \neq \emptyset$ and Neumann boundary conditions are imposed. Then, for any potential $q$ in $L_c^\infty (M) $, we have 
$$\lambda_1 (  q) \leq \lambda_1 (  c)=c,$$
where the equality holds if and only if $q=c$. Moreover, the constant potential $c$ is the only critical one of the functional $\lambda_1$ over $L_c^\infty (M)$. 

\item Assume that $\partial M \neq \emptyset$ and that Zero Dirichlet boundary conditions are imposed. Then the functional $\lambda_1$ does not admit any critical potential in $L_c^\infty (M)$.

\end{enumerate}
\end{theorem}

Higher eigenvalues are continuous but not differentiable in general. Nevertheless, perturbation theory enables us to prove that, for any function $u\in L^\infty (M)$, the function $t \mapsto \lambda_i (  q + tu)$ admits left and right derivatives at $t=0$ (see section 2.2). A generalized notion of criticality can be naturally defined as follows :
 
\begin {definition} A potential $q$ is said to be critical for the functional $ \lambda_i $ if, for any
$u \in L_*^\infty (M)$, the left and right derivatives of $t \mapsto \lambda_i (  q + tu)$ at $t=0$ have opposite signs, that is 
$${d\over {dt}} \lambda_i ( q + tu)\Big|_{t=0^+}\times {d\over {dt}} \lambda_i (  q + tu)\Big|_{t=0^-} \le 0.$$

\end{definition} 

It is immediate to check that $q$ is critical for $\lambda_i$ if and only if, for any $u \in L_*^\infty (M)$, one of the two following inequalities holds :
$$\lambda_i (  q + tu) \leq \lambda_i (  q) + o(t)\quad \mbox{as}\; t\to 0$$ 
or
$$\lambda_i (  q + tu) \geq \lambda_i (  q) + o(t)\quad \mbox{as}\; t\to 0.$$

In all the sequel, we will denote by $E_i (  q)$ the eigenspace corresponding to the $i$-th eigenvalue $\lambda_i (q)$ whose dimension coincides with the number of indices $j\in {\mathbb N}$ such that $ \lambda_j (  q)=\lambda_i (  q)$.

\medskip

As for the first eigenvalue, the functionals $ \lambda_i $, $ i\ge2$, admit no critical potentials under Dirichlet boundary conditions.  
\begin{theorem}
 Assume that $\partial M \neq \emptyset$ and that Zero Dirichlet boundary conditions are imposed. Then, $\forall i\in {\mathbb N}^ *$, the functional $\lambda_i$ does not admit any critical potential in $L_c^\infty (M)$.

\end{theorem}
Under the two remaining boundary conditions, the following theorem gives a necessary condition for a potential $q$ to be critical for the functional $ \lambda_i $. This condition is also sufficient for the indices $i$ such that $\lambda_i(  q)>\lambda_{i-1}(  q)$ or $\lambda_i (  q)<\lambda_{i + 1}(  q)$, which means that $\lambda_i(  q)$ is the first one or the last one in a cluster of equal eigenvalues.  

\begin{theorem}
Assume that either $\partial M = \emptyset$ or $\partial M \neq \emptyset$ and Neumann boundary conditions are imposed. Let $i$ be a positive integer.

If $q\in L_c^\infty (M)$ is a critical potential of the functional $\lambda_i$, then $q$ is smooth and
there exists a finite family of eigenfunctions $f_1 ,\ldots,f_k$ in $E_i (  q)$ such that $\sum_{1\leq j \leq k} f_j^2 = 1$.

Reciprocally, if $\lambda_i(  q)>\lambda_{i-1}(  q)$ or
$\lambda_i ( q)<\lambda_{i + 1}( q)$, and if there exists a family of eigenfunctions $f_1 ,\ldots,f_k \in E_i (  q)$ such that $\sum_{1\leq j \leq k} f_j^2 = 1$, then $q$ is a critical potential of the functional $\lambda_i$.

\end{theorem}
Note that the identity $\sum_{1\leq j \leq k} f_j^2 = 1$, with $f_1 ,\ldots,f_k\in E_i ( q)$, immediately implies another one (that we obtain from $ \Delta \sum_{1\leq j \leq k} f_j^2=0$):
$$q= \lambda_i( q)-\sum_{1\leq j \leq k} |\nabla f_j|^2,$$
from which we can deduce the smoothness of $q$.
\begin{remark}
1. The identity $\sum_{1\leq j \leq k} f_j^2 = 1$ with $-\Delta f_j +q f_j=\lambda_i(  q)f_j$, means that the map $f= (f_1, \ldots, f_k)$ from $M$ to the Euclidean sphere ${\mathbb S}^{k-1}$ is harmonic with energy density $|\nabla f|^2 = \lambda_i(  q) - q$ (see \cite{EL}). Hence, a necessary (and sometime sufficient) condition for a potential $q$ to be critical for the functional $\lambda_i$ is that the function $
\lambda_i(  q) - q$ is the energy density of a harmonic map from $M$ to a Euclidean sphere.

2. If one replaces the constraint on the mean value ${1\over V(M)}\int_Mqdv = c$ by the general constraint $\int_MF(q)dv = c,$ then the necessary and sufficient condition $\sum_{1\leq j \leq k} f_j^2 = 1$ of Theorem 1.3  becomes (even under Dirichlet boundary conditions) $\sum_{1\leq j \leq k} f_j^2 = F'(q)$. In particular, $q$ is a critical potential of the functional $\lambda_1$ if and only if $F'(q)\ge 0$ and $F'(q)^{1\over 2}$ is a first eigenfunction of $-\Delta+q$,
see \cite{AH, Ha} for a discussion of the case $F(q) = |q|^\alpha$.  
\end{remark}

Under each one of the boundary conditions we consider, a constant function can never be an eigenfunction associated to an eigenvalue $\lambda_i ( q)$ with $i\ge 2$. Hence, an immediate consequence of  Theorem 1.3 is the following
\begin{corollary}
If $q\in L_c^\infty (M)$ is a critical potential of the functional $\lambda_i$ with $i\ge 2$, then the eigenvalue $\lambda_i (  q)$ is degenerate, that is $\lambda_i (  q)=\lambda_{i-1} (  q)$ or $\lambda_i (  q)=\lambda_{i+1} (  q)$
\end{corollary}

If $\{f_1 ,\ldots,f_k\}$ is an $L^2$-orthonormal basis of $E_i (-\Delta )$, then the function $\sum_{1\leq j \leq k} f_j^2 $ is invariant under the isometry group of $M$. Indeed, for any isometry $\rho$ of $M$, $\{f_1\circ\rho ,\ldots,f_k\circ\rho\}$ is also an orthonomal basis of $E_i (-\Delta )$ and then, there exists a matrix $A\in O(d)$ such that $(f_1\circ\rho, \ldots, f_d\circ\rho)=A.(f_1, \ldots, f_d)$. In particular, if $M$ is homogeneous, that is, the isometry group acts transitively on $M$, then $\sum_{1\leq j \leq k} f_j^2$  would be constant. Another consequence of Theorem 1.3 is then the following

\begin{corollary}
 If $M$ is homogeneous, then constant potentials are critical
for all the functionals $\lambda_i$ such that $\lambda_i (-\Delta) < \lambda_{i+1}
(-\Delta)$ or $\lambda_i (-\Delta)> \lambda_{i-1}  (-\Delta)$.

\end{corollary}

Recall that Euclidean spheres, projective spaces and flat tori are examples of
homogeneous Riemannian spaces.

A potential $q \in L_c^\infty (M)$ is said to be a {\it local minimizer} (resp. {\it local maximizer}) of the functional $\lambda_i$ (in a weak sense) if, for any $u\in L_*^\infty (M)$, the function $t \mapsto \lambda_i (  q+t u)$ admits a local minimum (resp. maximum) at $t=0$. The result of Corollary 1.1 takes the following more precise form in the case of a local minimizer or maximizer.

\begin{theorem}
Let $q \in L_c^\infty (M)$ and $i\ge 2$.
\begin{enumerate} 
\item If $q$ is a local minimizer of the functional $\lambda_i$, then $\lambda_i(  q)=\lambda_{i-1}(  q)$.  
\item If $q$ is a local maximizer of the functional $\lambda_i$, then $\lambda_i(  q)=\lambda_{i+1}(  q)$. 
\end{enumerate}
\end{theorem}

Since the first eigenvalue is simple, we always have $\lambda_2( q)>\lambda_1( q)$. The previous results, applied to the functional $\lambda_2$ can be summarized as follows.

\begin{corollary}
Assume that either $\partial M = \emptyset$ or $\partial M \neq \emptyset$ and Neumann boundary conditions are imposed. 
A potential $q \in L_c^\infty (M)$ is critical for the functional $\lambda_2$ if and only if, $q$ is smooth, $\lambda_2( q)=\lambda_3( q)$ and there exist eigenfunctions $f_1 ,\ldots,f_k$ in $E_2 (  q)$ such that $\sum_{1\leq j \leq k} f_j^2 = 1$. 

Moreover, the functional $\lambda_2$ admits no local minimizers in  $L_c^\infty (M)$.
\end{corollary}

In \cite{EI}, Ilias and the first author have proved that, under some hypotheses on $M$, satisfied in particular by compact rank-one symmetric spaces, irreducible homogeneous Riemannian spaces and some flat tori, the constant potential $c$ is a global maximizer of $\lambda_2$ over $L_c^\infty (M)$. In \cite{EI2, EI3}, they studied the critical points of $\lambda_i$ considered as a functional on the set of Riemannian metrics of fixed volume on $M$.


\subsection{Critical potentials of the eigenvalue gaps functionals} $\quad$

\medskip

\noindent We consider now the eigenvalue gaps functionals $q\mapsto G_{ij}(q)=\lambda_j (  q) -\lambda_i(  q)$, where $ i$ and $ j$ are two distinct positive integers, and define their critical potentials as in Definition 1.1. These functionals are invariant under translations, that is $G_{ij}(q+c)=G_{ij}(q)$. Therefore, critical potentials of $G_{ij}$ with respect to fixed mean value deformations are also critical  with respect to arbitrary deformations. 

\begin{theorem}

If $q \in L_c^\infty (M)$ is a critical potential of the gap functional $G_{ij}=\lambda_j -\lambda_i$, then there exist a finite family of eigenfunctions $f_1 ,\ldots,f_k$ in $E_i(  q)$ and a finite family of eigenfunctions $g_1 ,\ldots,g_l$ in $E_{j}( q)$, such that $\sum_{1\leq p \leq k} f_p^2=\sum_{1\leq p \leq l} g_p^2 $.

Reciprocally, if $\lambda_i( q)<\lambda_{i+1}( q)$ and $\lambda_j( q)>\lambda_{j-1}( q)$, and if there exist $f_1 ,\ldots,f_k$ in $E_i(  q)$ and $g_1 ,\ldots,g_l$ in $E_{j}( q)$ such that $\sum_{1\leq p \leq k} f_p^2 = \sum_{1\leq p \leq l} g_p^2 $, then $q$ is a critical potential of $G_{ij}$.

\end{theorem}

In the particular case of the gap between two consecutive eigenvalues, we have the following 

\begin{corollary}
A potential $q \in L_c^\infty (M)$ is critical for the gap functional $G_{i,i+1}=\lambda_{i+1} -\lambda_i$ if and only if, either $\lambda_{i+1}( q)=\lambda_i( q)$, or there exist a family of eigenfunctions $f_1 ,\ldots,f_k$ in $E_i(  q)$ and a family of eigenfunctions $g_1 ,\ldots,g_l$ in $E_{i+1}( q)$, such that $\sum_{1\leq p \leq k} f_p^2=\sum_{1\leq p \leq l} g_p^2 $.
\end{corollary}

\begin{remark}
The characterization of critical potentials of $G_{ij}$ given in Theorem 1.5 remains valid under the constraint $\int_MF(q)dv=c$. 
\end{remark}

An immediate consequence of Theorem 1.5 is the following
\begin{corollary}
Let $q \in L_c^\infty (M)$ be a critical potential of the gap functional $G_{ij}=\lambda_j -\lambda_i$. If $\lambda_i( q)$ (resp. $\lambda_j (q)$) is nondegenerate, then $\lambda_j( q)$ (resp. $\lambda_i( q)$) is degenerate.

 \end{corollary}

The following is an immediate consequence of the discussion above concerning homogeneous Riemannian manifolds.
\begin{corollary}
If M is a homogeneous Riemannian manifold, then, for any positive integer $i$, constant  potentials are critical points of the gap functional $G_{i,i+1}=\lambda_{i+1} -\lambda_i$.
 \end{corollary}

Potentials $q$ such that $\lambda_{i+1}(  q)=\lambda_i( q)$ are of course global minimizers of the gap functional $G_{i,i+1}$. These potentials are also the only local minimizers of  $G_{i,i+1}$. Indeed, we have the following
\begin{theorem} 
 If $q \in L_c^\infty (M)$ is a local minimizer of the gap functional $G_{ij}=\lambda_j -\lambda_i$, then, either $\lambda_i( q)=\lambda_{i+1}( q)$, or $\lambda_j( q)=\lambda_{j-1}( q)$.
If $q$ is a local maximizer of $G_{ij}$, then, either $\lambda_i( q)=\lambda_{i-1}( q)$, or $\lambda_j( q)=\lambda_{j+1}( q)$.

In particular, $q $ is a  local minimizer of the gap functional $G_{i,i+1}=\lambda_{i+1} -\lambda_i$ if and only if $G_{i,i+1}(q)=0$.
\end{theorem}

Finally, let us apply the results of this section to the first gap $G_{1,2}$.

\begin{corollary}
 A potential $q \in L_c^\infty (M)$ is critical for the gap functional $G_{1,2}=\lambda_{2} -\lambda_1$ if and only if $\lambda_{2}( q)$ is degenerate and there exists a family of eigenfunctions $g_1 ,\ldots,g_l$ in $E_{2}( q)$ such that $\sum_{1\leq j \leq l} g_j^2 =f^2$, where $f$ is a basis of $E_{1}( q)$.

The functional $G_{1,2}$ does not admit any local minimizer in $L_c^\infty (M)$.

 \end{corollary}

The authors wish to thank the referee for his valuable remarks.


\section {Proof of Results}\label{2}
\subsection{Variation Formula and proof of Theorem 1.1}
Given on $M$ a potential $q$ and a function $u\in L^\infty(M)$, we consider  the family of operators $-\Delta + q +tu$. Suppose that $\Lambda (t)$ is a differentiable family of eigenvalues of $-\Delta + q +tu$ and that $f_t$ is a differentiable family of corresponding normalized eigenfunctions, that is, $\forall t$,
$$(-\Delta + q +tu) f_t=\Lambda (t) f_t,$$
and
$$\int_{M}f_t^2dv =1,$$
with $f_t\big|_{\partial M}=0$ or ${\partial f_t\over \partial \nu}\big|_{\partial M}=0$ if $\partial M\neq\emptyset$.
The following formula, giving the derivative of $\Lambda$, is already known at least in the case of Euclidean domains with Dirichlet boundary conditions. 

\begin {proposition} $${\Lambda}'(0)=\int_M u f_0^2 dv.$$

\end {proposition}

\begin{proof}
 First, we have, for all $t$,
$$\Lambda (t)=\Lambda (t)\int_{M} (f_t)^2 dv = \int_M f_t(-\Delta + q + tu)f_t \, dv.$$
Differentiating at $t=0$, we get
$$\Lambda'(0) = {d\over dt} \biggl(\int_M f_t (-\Delta + q)f_t \, dv + 
t\int_M u(f_t)^2 dv \biggr)\Big|_{t=0}. $$
Now, noticing that the function ${d\over dt}f_t\big|_{t=0}$ satisfies the same boundary conditions as $f_0$ in case $\partial M \neq \emptyset$, and using integration by parts, we obtain  

\begin{eqnarray}
\nonumber{}{d\over dt} \int_M f_t (-\Delta + q)f_t \,dv \Big|_{t=0}&=& 2 \int_{M}(-\Delta + q)f_0\, {d\over dt}f_t\Big|_{t=0} dv \\
\nonumber{} &=& 2\Lambda (0) \int_{M}f_0\;{d\over dt}f_t\Big|_{t=0} dv \\ 
\nonumber{} &=& \Lambda (0)\;{d\over dt} \int_{M}f_t^2dv \Big|_{t=0} = 0.
\end{eqnarray}

On the other hand, we have
\begin{eqnarray} 
\nonumber{} {d\over dt} \Big( t\int_M uf_t^2 dv\Big)  \Big|_{t=0} &=& \int_M uf_0^2dv +\Big(t\int_M u{d\over dt}f_t^2 dv \Big) \Big|_{t=0} \\
\nonumber{}  &=& \int_M uf_0^2v_g.
\end{eqnarray}
Finally, ${\Lambda}' (0)=\int_M uf_0^2 dv$.
\end{proof}

\begin{proof}
(of Theorem 1.1.)
(i) First, let us show that constant potentials are maximizing for $\lambda_1$. Indeed, let $c$ be a constant potential and let $q$ be an arbitrary one in $L_c^\infty(M)$. From the variational characterization of $\lambda_1(-\Delta+q)$ in the case $\partial M=\emptyset$ as well as in the case of Neumann boundary conditions, we get
\begin{eqnarray}
\nonumber{} \lambda_1(-\Delta+q) &=& \inf_{f \in H^1(M)} {\int_M (|\nabla f|^2 +qf^2) dv\over \| f \|^2_{L^2(M)}} \\
 \nonumber{} & \leq &{\int_M (|\nabla 1|^2 +q1^2) dv\over \| 1 \|^2_{L^2(M)}}={\int_M q\, dv \over V(M)}=c.
\end{eqnarray}  
Hence, $\lambda_1 (  q)\leq \lambda_1 (c)$ and the constant potential $c$
maximizes the functional $\lambda_1$ on $L_c^\infty(M)$. In particular, constant potentials are critical for this  functional. 

Now, suppose that $q\in L_c^\infty(M)$ is a critical potential for $\lambda_1$. For any $u \in L_*^\infty(M)$,
we consider a differentiable family $f_t$ of normalized eigenfunctions corresponding to the first eigenvalue of $(-\Delta + q + tu)$ and apply the variation formula above to obtain: 
$${d\over dt}\lambda_1(  q + tu)\Big|_{t=0} = \int_{M} uf_0^2\, dv.$$
Hence, $\int_M uf_0^2 dv =0$ for any $u \in L_*^\infty (M)$, which implies that $f_0$ is constant on $M$. Since $(-\Delta + q)f_0 =q f_0 = \lambda_1(  q)f_0$, the potential $q$ must be constant on $ M$.

(ii) Let $f_0$ be the first nonnegative Dirichlet eigenfunction of $-\Delta + q$ satisfying $\int_M f_0^2\, dv=1$. The function $u=V(M) f_0^2 - 1$ belongs to $ L_*^\infty(M)$ and we have
$${d\over dt}\lambda_1(  q + tu)\Big|_{t=0} = \int_{M} uf_0^2\, dv=V(M) \int_M f_0^4\, dv-1>0,$$
where the last inequality comes from  Cauchy-Schwarz inequality and the fact that $f_0$ is not constant (recall that $f_0\big|_{\partial M}=0$). Therefore, the potential $q$ is not critical for $\lambda_1$.
\end{proof}


\subsection{Characterization of critical potentials}
Let $i$ be a positive integer and let $m\ge 1$ be the dimension of the eigenspace $E_i(  q)$ associated to the eigenvalue $ \lambda_i(  q)$. For any function $u\in  L_*^\infty(M)$, perturbation theory of unbounded self-adjoint operators (see for instance Kato's book \cite{K}) that we apply to the one parameter family of operators $-\Delta + q + tu$, tells us that, there exists a family of $m$ eigenfunctions
$f_{1,t},\ldots,f_{m,t}$ associated with a family of $m$ (non ordered) eigenvalues
$\Lambda_1(t),\ldots,\Lambda_m(t)$ of $-\Delta + q + tu$, all depending analytically in $t$ in some interval $(-\varepsilon, \varepsilon)$, and satisfying
\begin{itemize}
\item  $\Lambda_1(0)=\cdots =\Lambda_m(0)=\lambda_i(  q)$, 
\item  $\forall t\in (-\varepsilon, \varepsilon)$, the $m$ functions $f_{1,t},\ldots,f_{m,t}$ are 
orthonormal in $L^2(M)$. 
\end{itemize}
From this, one can easily deduce the existence of two integers $k\le m$ and $l\le m$, and a small $\delta>0$ such that
$$ \lambda_i(  q+tu) =\left\{ \begin{array}{l}
\Lambda_k(t)\ \hbox{if}\ t\in (-\delta, 0)\\
\\
\Lambda_l(t)\ \hbox{if}\ \ t \in (0,\delta ).
\end{array}\right.$$
Hence, the function $t\mapsto\lambda_i(  q+tu)$ admits a left sided and a right
sided derivatives at $t=0$ with
$${d\over dt} \lambda_i(  q+tu)\Big|_ {t=0^-}
=  \Lambda'_k (0)=\int_M uf_{k,0}^2 dv$$
and
$${d\over dt} \lambda_i(  q+tu) \Big|_ {t=0^+}
=  \Lambda'_l (0)=\int_M uf_{l,0}^2 dv.$$

To any function $u\in L^\infty_*(M)$ and any integer $i\in {\mathbb N} $, we associate the quadratic form $Q_u^i$ on $E_i(  q)$ defined by 
$$Q_u^i(f)=\int_M uf^2 dv.$$
The corresponding symmetric linear transformation $L_u^i: E_i(  q)\rightarrow E_i(  q)$ is given by
$$L_u^i (f)=P_i(uf),$$
where $P_i: L^2(M)\to E_i(  q)$ is the orthogonal projection of $L^2(M)$ onto $E_i(  q)$. 
 
It follows immediately that

\begin{proposition}
If the potential $q$ is critical for the functional $\lambda_i$, then,  $\forall u\in L^\infty_*(M)$, the quadratic form $Q_u^i(f)=\int_M uf^2 dv$ is indefinite on the eigenspace $E_i(  q)$.
\end{proposition}

The following lemma enables us to establish a converse to this proposition. 
\begin{lemma}
$\forall k, l \le m$, we have
$$\int_Mu f_{k,0}f_{l,0}\, dv =\left\{ \begin{array}{l}
0\ \hbox{if}\ \ k\neq l\\
\\
\Lambda'_k(0)\ \hbox{if}\ \ k=l.
\end{array}\right.$$
In other words, $\Lambda'_1(0),\ldots,\Lambda'_m(0)$ are the eigenvalues of the symmetric linear transformation $L_u^i: E_i(  q)\rightarrow E_i(  q)$ and the functions $f_{1,0},\ldots,f_{m,0}$ constitute an orthonormal eigenbasis of $L_u^i$.

\end{lemma}

\begin{proof} 
Differentiating at $t=0$ the equality $(-\Delta + q+tu)f_{k,t} = \Lambda_k(t) f_{k,t}$, we obtain
$$uf_{k,0} +(-\Delta + q)\frac{d}{dt}f_{k,t}\Big|_{t=0} = \Lambda'_k(0) f_{k,0}+\Lambda_k(0)\frac{d}{dt}f_{k,t}\Big|_{t=0},$$
and then,
\begin{eqnarray} 
\nonumber{} \int_Muf_{k,0}f_{l,0}\, dv = \Lambda'_k(0) \int_Mf_{k,0} f_{l,0}\, dv  &+& \Lambda_k(0)\int_Mf_{l,0}\frac{d}{dt}f_{k,t}\Big|_{t=0}dv\\
\nonumber{}  &-&\int_Mf_{l,0}(-\Delta + q)\frac{d}{dt}f_{k,t}\Big|_{t=0}dv. 
\end{eqnarray}
Integration by parts gives, after noticing that $\Lambda_k(0)=\Lambda_l(0)=\lambda_i(  q)$ and that the functions $\frac{d}{dt}f_{k,t}\big|_{t=0}$ satisfy the considered boundary conditions, 
\begin{eqnarray} 
\nonumber{} \int_Mf_{l,0}(-\Delta + q)\frac{d}{dt}f_{k,t}\Big|_{t=0}dv &=& \int_M\frac{d}{dt}f_{k,t}\Big|_{t=0}(-\Delta + q)f_{l,0}\, dv\\
\nonumber{}  &=& \Lambda_k(0)\int_Mf_{l,0}\frac{d}{dt}f_{k,t}\Big|_{t=0}dv,
\end{eqnarray}
and finally,
$$\int_Muf_{k,0}f_{l,0}\, dv= \Lambda'_k(0) \int_Mf_{k,0} f_{l,0}\, dv =\Lambda'_k(0)\delta_{kl}.$$
\end{proof}

\begin{proposition} Assume that $\lambda_i(  q)>\lambda_{i-1}(  q)$ or
$\lambda_i (  q)<\lambda_{i + 1}(  q)$.
Then the following conditions are equivalent:
\begin{itemize}
\item[i)]  the potential $q$ is critical for $\lambda_i$
\item[ii)]  $\forall u\in L^\infty_*(M)$, the quadratic form $Q_u^i(f)=\int_M uf^2 dv$ is indefinite on the eigenspace $E_i(  q)$.

\item[iii)]   $\forall u\in L^\infty_*(M)$, the linear transformation $L_u^i$ admits eigenvalues of both signs. 

\end{itemize}
\end{proposition}

\begin{proof} Conditions (ii) and (iii) are clearly equivalent and the fact that (i) implies (ii) was established in Proposition 2.2. Let us show that (iii) implies (i). Assume that $\lambda_i(  q)>\lambda_{i-1}(  q)$ and let $u\in L^\infty_*(M)$ and $\Lambda_1(t),\ldots,\Lambda_m(t)$ be as above. For small $t$, we will have, for continuity reasons, $\forall k\le m$, $\Lambda_k(t) >\lambda_{i-1}(  q+tu)$ and then, $ \lambda_i(  q+tu) \le \Lambda_k(t)$. Since $ \lambda_i(  q+tu) \in \{\Lambda_1(t),\ldots,\Lambda_m(t)\}$, we get 
$$ \lambda_i(  q+tu) = \min_{k\le m} \Lambda_k(t). $$
It follows that
$$\frac{d}{dt}\lambda_i(  q+tu)\Big| _{t=0^-} = \max_{k\le m} \Lambda'_k(0) $$
and
$$\frac{d}{dt}\lambda_i(  q+tu)\Big| _{t=0^+} = \min_{k\le m} \Lambda'_k(0) .$$
Thanks to Lemma 2.1, Condition (iii) implies that $\min_{k\le m} \Lambda'_k(0) \le 0\le 
 \max_{k\le m} \Lambda'_k(0) $ which implies the criticality of $q$. 

The case $\lambda_i (  q)<\lambda_{i + 1}(  q)$ can be treated in a similar manner.
\end{proof}


\subsection{Proof of Theorems 1.2 and 1.3} Let $q$ be a potential in $L^\infty_c(M)$. To prove Theorem 1.2 we first notice that, since $f\big|_{\partial M}=0$ for any $f\in E_i(  q)$, the constant function 1 does not belong to the vector space $F$ generated in $L^2(M)$ by $\{f^2\, ; f \in E_i(  q)\}$. Hence, there exists a function $u$ orthogonal to $F$ and such that $\langle u, 1\rangle _{L^2(M)} <0$. The function $u_0= u- \bar u$ belongs to $L^\infty_*(M)$ and the quadratic form $Q^i_{u_0}(f) =  \int_M u_0f^2 dv = - \bar u  \|f\|^2_{L^2(M)}$ is positive definite on $E_i(q)$. Hence, the potential $q$ is not critical for $\lambda_i$ (see Proposition 2.2).

The proof of Theorem 1.3 follows directly from the two propositions above and the following lemma.

\begin{lemma}
Let $i$ be a positive integer. The two following conditions are equivalent:
\begin{itemize}
\item [i)] $\forall u\in L^\infty_*(M)$, the quadratic form $Q_u^i(f)=\int_M uf^2 dv$ is indefinite on the eigenspace $E_i(  q)$.
\item [ii)]  there exists a family of eigenfunctions $f_1 ,\ldots,f_k$ in $E_i (  q)$ such that $\sum_{1\leq j \leq k} f_j^2 = 1$.
\end{itemize}
\end{lemma}

\begin{proof} To see that (i) implies (ii) we introduce the convex cone $C$ generated in $L^2(M)$ by the set $\{f^2\, ; f \in E_i(  q)\}$, that is $C =\{\sum_{j \in J}f_j^2\, ; f_j \in E_i( q),J \subset {\mathbb N} \, ,\ J\,
\hbox{ is finite} \}$. Condition (ii) is then equivalent to the fact that the constant function 1 belongs to $C$. Let us suppose, for a contradiction, that $1\notin C$. Then, applying classical separation theorems (in the finite dimensional vector subspace of $L^2(M)$ generated by $\{f^2\, ; f \in E_i(  q)\}$ and $1$, see \cite{R}), we prove the existence of a function
$u \in L^2(M)$ such that $\bar u= \frac{1}{V(M)} \int_M u\cdot 1\, dv < 0$ and $\int_M uf^2 dv \ge 0$ for any $f \in C$.
Hence, the function $u_0= u- \bar u$ belongs to $L^\infty_*(M)$ and satisfies, $\forall f \in E_i(q)$, 
$$Q^i_{u_0}(f) =  \int_M uf^2dv -  \frac{1}{V(M)} \int_M u \, dv \int_M f^2dv \ge - \bar u  \|f\|^2_{L^2(M)}.$$
The quadratic form $Q^i_{u_0}$ is then positive definite which contradicts (i) (see Proposition 2.2). 

Reciprocally, the existence of $f_1 ,\ldots,f_k$ in $E_i (  q)$ satisfying $\sum_{1\leq j \leq k} f_j^2 = 1$ implies that, $\forall u\in L^\infty_*(M)$,
$$\sum_{j\le k} Q^i_u(f_j) = \sum_{j\le k} \int_M uf_j^2 dv =\int_M u =0,$$
which implies that the quadratic form $Q^i_u$ is indefinite on $E_i(  q)$.
\end{proof}

Finally, let us check that the condition $\sum_{1\leq j \leq k} f_j^2 = 1$, with $f_j \in E_i(  q)$, implies that $q$ is smooth. Indeed, since $q\in L^\infty (M)$, we have, for any eigenfunction $f \in E_i(  q)$, $\Delta f \in L^2(M)$ and then, $f\in H^{2,2}(M)$. Using standard regularity theory and Sobolev embeddings (see, for instance, \cite{J}), we obtain by an elementary iteration, that  $f\in H^{2,p}(M)$ for some $p>n$, and, then, $f\in C^1(M)$. From $\sum_{1\leq j \leq k} f_j^2 = 1$ and $\Delta \sum_{1\leq j \leq k} f_j^2 = 0$, we get 
$$q= \lambda_i( q)-\sum_{1\leq j \leq k} |\nabla f_j|^2,$$
which implies that $q$ is continuous. Again, elliptic regularity theory tells us that the eigenfunctions of $-\Delta + q$ are actually smooth, and, hence, $q$ is smooth.


\subsection{Proof of Theorem 1.4}

Assume that the potential $q$ is a local minimizer of the functional $\lambda_{i}$ on $L^\infty_c(M)$ and let us suppose for a contradiction that $\lambda_i(  q)>\lambda_{i-1}(  q)$. Let $u$ be a function in $L^\infty_*(M)$ and let $\Lambda_1(t),\ldots,\Lambda_m(t)$ be a family of $m$ eigenvalues of $-\Delta + q + tu$, where $m$ is the multiplicity of $\lambda_i(  q)$, depending analytically in $t$ and such that $\Lambda_1(0)=\cdots=\Lambda_m(0)=\lambda_i(  q)$. For continuity reasons, we have, for sufficiently small $t$ and any $k\le m$, $\Lambda_k(t)>\lambda_{i-1}(  q +tu)$. Hence, $\forall k\le m$ and $\forall t$ sufficiently small, 
$$\Lambda_k(t)\ge\lambda_{i}(  q +tu)\ge\lambda_{i}(  q) =\Lambda_k(0).$$  
Consequently, $\forall k\le m$, $\Lambda'_k(0)=0$. Applying Lemma 2.1 above we deduce that the symmetric linear transformation $L^i_u$ and then the quadratic form $Q^i_u$ is identically zero on the eigenspace $E_i (  q)$. Therefore, $\forall u\in L^\infty_*(M)$ and $\forall f \in E_i (  q)$, we have $\int_M uf^2v_g=0$. In conclusion, $\forall f \in E_i (  q)$, $f$ is constant on $M$ which is impossible for $i\ge 2$. 
The same arguments work to prove Assertion (ii).


\subsection{Proof of Theorem 1.5}

Let $q$ be a potential and let $i$ and $j$ be two distinct positive integers such that $\lambda_i( q)\neq\lambda_j( q)$. We denote by $m$ (resp. $n$) the dimension of the eigenspace $E_i(  q)$ (resp. $E_j(  q)$). Given a function $u$ in $L_*^\infty (M)$, we consider, as above, $m$ (resp. $n$) $L^2(M)$-orthonormal families of eigenfunctions $f_{1,t},\ldots,f_{m,t}$ (resp. $g_{1,t},\ldots,g_{n,t}$) associated with $m$ (resp. $n$) families of eigenvalues $\Lambda_1(t),\ldots,\Lambda_m(t)$ (resp. $\Gamma_1(t),\ldots,\Gamma_n(t)$) of $ -\Delta+ q + tu$, all depending analytically in $t\in (-\varepsilon, \varepsilon)$, such that $\Lambda_1(0)=\cdots=\Lambda_m(0)=\lambda_i(  q)$ (resp. $\Gamma_1(0)=\cdots=\Gamma_n(0)=\lambda_j( q)$). Hence, there exist four integers $k\le m$, $k'\le m$, $l\le n$ and $l'\le n$, such that
\begin{eqnarray}
\nonumber{} {d\over dt} (\lambda_j-\lambda_i)( q+tu)\Big|_ {t=0^-}&=&\Gamma'_l(0)-\Lambda'_k (0)\\
\nonumber{}&=&\int_M u( g_{l,0}^2 -f_{k,0}^2)dv
\end{eqnarray}
and
\begin{eqnarray}
\nonumber{} {d\over dt} (\lambda_j-\lambda_i)(  q+tu)\Big|_ {t=0^+}&=&\Gamma'_{l'}(0)-\Lambda'_{k'} (0)\\
\nonumber{}&=&\int_M u( g_{l',0}^2 -f_{k',0}^2)dv.
\end{eqnarray}
Recall that (Lemma 2.1) the eigenfunctions $f_{1,0},\ldots,f_{m,0}$ (resp. $g_{1,0},\ldots,g_{n,0}$) constitutes an $L^2(M)$-orthonormal basis of $E_i(  q)$ (resp. $E_j(  q)$) which diagonalizes the quadratic form $Q_u^i$ (resp. $Q_u^j$). Therefore, the family $(f_{k,0} \otimes g_{l,0})_{k\le m\, ,\, l\le n}$ constitutes a basis of the space $E_i(  q)\otimes E_j(  q)$ which diagonalizes the quadratic form $S_u^{i,j}$ given by 
\begin{eqnarray}
\nonumber{} S_u^{i,j}(f\otimes g)&=&\|f\|^2_{L^2(M)} Q_u^j(g) - \|g\|^2_{L^2(M)} Q_u^i(f) \\
\nonumber{} &=&  \int_M u( \|f\|^2_{L^2(M)}g^2- \|g\|^2_{L^2(M)}f^2)dv.
\end{eqnarray}
The corresponding eigenvalues are $ (\Gamma'_l(0)-\Lambda'_k (0))_{k\le m\, ,\, l\le n}$. The criticality of $q$ for $\lambda_j - \lambda_i$ then implies that this quadratic form admits eigenvalues of both signs, which means that it is indefinite.

On the other hand, in the case where $\lambda_i( q)<\lambda_{i+1}( q)$ and $\lambda_j( q)>\lambda_{j-1}( q)$, we have, as in the proof of Proposition 2.3,  for sufficiently small $t$,
$\lambda_i( q+tu) = \max_{k\le m}\Lambda_k(t)$ and $\lambda_j( q+tu) =\min _{l\le n}\Gamma_l(t)$, which yields
\begin{eqnarray}
\nonumber{} {d\over dt} (\lambda_j-\lambda_i)(  q+tu)\Big|_ {t=0^-}&=&\max _{l\le n}\Gamma'_l(0)-\min_{k\le m}\Lambda'_k(0)\\
\nonumber{}&=&\max_{k\le m\, ,\, l\le n}(\Gamma'_l(0)-\Lambda'_k(0))
\end{eqnarray}
and 
\begin{eqnarray}
\nonumber{} {d\over dt} (\lambda_j-\lambda_i)(  q+tu)\Big|_ {t=0^+}&=&\min _{l\le n}\Gamma'_l(0)-\max_{k\le m}\Lambda'_k(0)\\
\nonumber{}&=&\min_{k\le m\, ,\, l\le n}(\Gamma'_l(0)-\Lambda'_k(0)).
\end{eqnarray}

One deduces the following

\begin{proposition}
If the potential $q\in L_c^\infty (M)$ is critical for the functional $G_{ij}=\lambda_{j}-\lambda_{i}$, then,  $\forall u\in L_*^\infty (M)$, the quadratic form $S_u^{i,j}$ is indefinite on $E_i(  q)\otimes E_j(  q)$.

Reciprocally, if $\lambda_i( q)<\lambda_{i+1}( q)$ and $\lambda_j( q)>\lambda_{j-1}( q)$, and if, $\forall u\in L_*^\infty (M)$,  the quadratic form $S_u^{i,j}(g)$ is indefinite on $E_i(  q)\otimes E_j(  q)$, then $q$ is a critical potential of the functional $G_{ij}$.
\end{proposition}

The following lemma will completes the proof of Theorem 1.5
\begin{lemma}
The two following conditions are equivalent:
\begin{itemize}
\item [i)] $\forall u\in L_*^\infty (M)$, the quadratic form $S_u^{i,j}$ is indefinite on $E_i(  q)\otimes E_j(  q)$.
\item [ii)]  there exist a finite family of eigenfunctions $f_1 ,\ldots,f_k$ in $E_i(  q)$ and a finite family of eigenfunctions $g_1 ,\ldots,g_l$ in $E_{j}( q)$, such that $\sum_{1\leq p \leq k} f_p^2=\sum_{1\leq p \leq l} g_p^2 $.
\end{itemize}
\end{lemma}
The proof of this lemma is similar to that of Lemma 2.2. Here, we consider the two convex cones $C_i $ and  $C_j$ in $L^2 (M)$ generated respectively by $\left\{f^2\, ;\, f \in E_i(q)\, ,f\neq 0\right\}$ and $\left\{g^2\, ;\, g \in E_j(q)\, ,g\neq 0\right\}$. Condition (ii) is then equivalent to the fact that these two cones admit a nontrivial intersection. As in the proof of Lemma 2.2, separation theorems enable us to prove that, if $C_i\cap C_j=\emptyset$, then there exists a function $u$ such that $\int_M u f^2 dv <0$ for any $f \in E_i(q)$, and $\int_M u g^2 dv \ge 0$ for any $f \in E_j(q)$, which implies that  $S_u^{i,j}$ is positive definite on $E_i(  q)\otimes E_j(  q)$. Since $S_1^{i,j} = 0$, we have,   $S_u^{i,j}=S_{u_0}^{i,j}$ with $u_0= u- \bar u \in L_*^\infty (M)$. Proposition 2.4 enables us to conclude. 

Reciprocally, assume the existence of $f_1 ,\ldots,f_k \in E_i(  q)$ and $g_1 ,\ldots,g_l \in E_{j}( q)$ satisfying $\sum_{1\leq p \leq k} f_p^2=\sum_{1\leq p \leq l} g_p^2 $. Then, $\forall u\in L_*^\infty (M)$, 
$$\sum_{1\leq p \leq k} \, \sum_{1\leq p' \leq l }S_u^{i,j}(f_p\otimes g_{p'}) =\cdots =0,$$
which implies that $S_u^{i,j}$ is indefinite on $E_i(  q)\otimes E_j(  q)$.

\subsection{Proof of Theorem 1.6}
Let $q$ be a local minimizer of $G_{ij}=\lambda_{j}-\lambda_{i}$ and let us suppose for a contradiction that $\lambda_i( q)<\lambda_{i+1}(  q)$ and $\lambda_j(  q)>\lambda_{j-1}(  q)$. Given a function $u$ in $L_*^\infty (M)$, we consider, as above, $m$ (resp. $n$) families of eigenvalues $\Lambda_1(t),\ldots,\Lambda_m(t)$ (resp. $\Gamma_1(t),\ldots,\Gamma_n(t)$) of $ -\Delta+ q + tu$, with $m=\dim E_i(q)$ and $n=\dim E_j(q)$, such that $\Lambda_1(0)=\cdots=\Lambda_m(0)=\lambda_i(  q)$ and $\Gamma_1(0)=\cdots=\Gamma_n(0)=\lambda_j( q)$. As in the proof of Theorem 1.4, we will have for sufficiently small $t$, $\lambda_i( q+tu) = \max_{k\le m} \Lambda_k(t)$ and $\lambda_j( q+tu) = \min_{l\le n} \Gamma_l(t)$. Hence, $\forall k\le m$ and $l\le n$,
\begin{eqnarray}
\nonumber{} \Gamma_l(t)- \Lambda_k(t)&\ge& \lambda_j( q+tu) -\lambda_i( q+tu) =G_{ij}(q+tu)\\
\nonumber{}&\ge& G_{ij}(q) =\Gamma_l(0)- \Lambda_k(0).
\end{eqnarray}
It follows that, $\forall k\le m$ and $l\le n$, $\Gamma'_l(0)- \Lambda'_k(0)=0$ and, then, the quadratic form $S_u^{i,j}$ is identically zero on $E_i(  q)\otimes E_j(  q)$ (recall that $\Gamma'_l(0)- \Lambda'_k(0)$ are the eigenvalues of $S_u^{i,j}$). This implies that, $\forall f\in E_i(  q)$ and $\forall g\in E_j(  q)$, the function $ \|f\|^2_{L^2(M)}g^2- \|g\|^2_{L^2(M)}f^2$ is constant equal to zero (since its integral vanishes) which is clearly impossible unless $i=j$.

\end{document}